\documentclass[titlepage,twoside,12pt]{article}
\usepackage{amssymb}
\usepackage{amsfonts}
\begin{document}

{\bf \Large Comments on a recent item by Yu I Manin} \\ \\

{\bf Elem\'{e}r E Rosinger} \\ \\
Department of Mathematics \\
and Applied Mathematics \\
University of Pretoria \\
Pretoria \\
0002 South Africa \\
eerosinger@hotmail.com \\ \\

{\bf Abstract} \\

The comments relate to the often overlooked, if not in fact intentionally disregarded depths
of what the so called internal aspects of mathematical knowledge may involve, depths
concerning among others issues such as its unreasonable effectiveness in natural sciences, to
use the terms of Eugene Wigner, suggested more than four decades ago. \\ \\

{\bf 1. On what Appears to be a Modern Trend ...} \\

It seems appropriate to start with the following citation from the European-American
philosopher of history, Eric Voegelin (1901-1985), who in the first chapter "The emergence of
secularized history : Bossuet and Voltaire" of his book "From Enlightenment to Revolution",
writes referring to an essay of Voltaire :

\begin{quote}

"... the Essay implies the idea that universality in historiography can be achieved by
completeness, and insofar as the Essay implies this identification it rather opens than solves
the problem of universality. By completeness one can achieve an encyclopedia, but not
automatically a unity of meaning. \\
... it reveals the weakness of all universal histories ... : the impossibility of finding a
meaning ... \\
The meaning, of course, cannot be found ... \\
The ideal of empirical completeness, however, cannot be more than a transitory position in the
movement of ideas. As soon as the question is raised why one should know with any degree of
completeness whatever has happened in the existence of mankind in time, the curiosity shop is
revealed as senseless. Encyclopedic knowledge, collected in handbooks, has to be moved into
the functional position of a collection of materials which ultimately might become of
importance for a relevant interpretation of history ..."

\end{quote}

Beyond, and far more important than any partial or complete encyclopedia - or "curiosity shop",
in terms of Voegelin - appear to be concerns about Mathematics, such as that expressed related
to its apparent {\it unreasonable effectiveness} in the Natural Sciences, [4]. \\
And in this regard, any wealth of classifications, divisions, or for that matter, endless
subdivisions of Mathematics, among them in Algebra or Geometry, in left brain versus right
brain activities, and so on and so on, has a focus that may recall a spectator's view point,
one who is sophisticated enough to discriminate with elaborate and on occasion astute comments,
yet not interested in, even if possibly capable of, following further and further "up stream",
towards the sources of mathematical thinking. Those very sources which, among others, make it
so unreasonably effective ... \\

In this regard, as it happens, there is no mention in [2] of that major issue brought up in
[4]. Instead, and as far as the deeper realms underlying mathematics thinking may be concerned,
mention is made of the extraordinary abilities of Ramanujan in dealing with "infinite
identities", abilities which remained totally mysterious even to G Hardy, a distinguished
mathematician at the time, and one who had the occasion to know and work with Ramanujan more
closely for several years. \\ \\

{\bf 2. Where and How Does It Happen ?} \\

It is typical of modern scientists, and not only mathematicians, that when involved in
research, or possibly even in the history or philosophy of science, they consider the fact or
event of their own thinking as totally outside of absolutely everything. So completely outside,
that it is simply not in this, or for that matter, not in any other given Universe ... \\
Of course, various psychological reflections may on occasion be made about the scientific
thinking process. And the feed-back effect upon one's thinking of the ongoing research
activity is readily acknowledged, actually sought out as most valuable, and seen in fact as a
highly necessary learning process. \\
Yet no one seems to bother much nowadays about the deeper sources of one's thinking. Instead,
one tends simply to consider oneself gifted and/or lucky, plus of course, as benefiting from
one's long hard work. \\

And with that, the issue of the fact or event of one's thinking tends to be closed for
evermore ... \\

Or is it ? \\

Not very much long ago, that is, less than four centuries earlier, and at the dawn of what we
tend to consider as the beginnings of modern science, Descartes believed to have pinpointed
the answer to this rather fundamental issue. He saw Creation divided up into two very much
distinct and disjoint realms. Namely, they were "res cogitans" and "res extensa", the former
where all thinking happens, and the latter the physical world accessible to our five usual
senses. \\
Of course, like nearly everybody else in his time, Descartes was deeply religious. Therefore,
according to his scheme, Creation may have been so sharply divided, yet underlying that
division was in his view the Creator himself, in all of his uniqueness and unity. And then,
res cogitans and res extensa were in a self-evident manner but the two branches growing out of
one single, all supporting, and unifying tree trunk ... \\

As Western scientists kept abandoning more and more the religious views of a Descartes, the
two Cartesian realms got severed from their single and unifying trunk, and thus they led to
the difficult and still not solved "mind-body" problem, a problem that is bound to include
that long untractable dichotomy which a Descartes did not have to face, and thus, deal
with ... \\

Strangely enough, for modern scientists the typical way to overcome that mind-body dichotomy
in their day to day research activity is simply by a decision, often taken by default and
less than consciously. Namely, they simply act as if saying to themselves :

\begin{quote}

Never mind about mind ... !

\end{quote}

And yet, no matter how widespread and entrenched such an approach may be, there are some most
simple questions, especially in modern Physics, which point to the grave inadequacy of such a
forgetful approach to the Cartesian twin realms of "res cogitans" and "res extensa". Indeed,
as mentioned for instance in [3], one can obviously ask questions such the following ones. \\

{\bf 2.1. Within Newtonian Mechanics.} Instant action at arbitrary distance, such as in the
case of gravitation, is one of the basic assumptions of Newtonian Mechanics. This does not
appear to conflict with the fact that we can think instantly and simultaneously about
phenomena no matter how far apart from one another in space and/or in time. However, absolute
space is also a basic assumption of Newtonian mechanics. And it is supposed to contain
absolutely everything that may exist in Creation, be it in the past, present or future.
Consequently, it is supposed to contain, among others, the physical body of the thinking
scientist as well. Yet it is not equally clear whether it also contains scientific thinking
itself which, traditionally, is assumed to be totally outside and independent of all phenomena
under its consideration, therefore in particular, totally outside and independent of the
Newtonian absolute space, and perhaps also of absolute time. \\

And then the question arises : where and how does such a scientific thinking take place or
happen ? \\

{\bf 2.2. Within Einstein's Mechanics.} In Special and General Relativity a basic assumption
is that there cannot be any propagation of action faster than light. Yet just like in the case
we happen to think in terms of Newtonian Mechanics, our thinking in terms of Einstein's
Mechanics can again instantly and simultaneously be about phenomena no matter how far apart
from one another in space and/or time. \\

Consequently, the question arises : given the mentioned relativistic limitation, how and where
does such a thinking happen ? \\

{\bf 2.3. With Quantum Mechanics.} Let us consider the classical EPR entanglement phenomenon,
and for simplicity, do so in terms of quantum computation. For that purpose it suffices to
consider double qubits, that is, elements of $\mathbb{C}^2 \bigotimes \mathbb{C}^2$, such as
for instance the pair \\

$~~~~~~~~~~~~~~~ |~ \omega_{00} > ~=~ |~ 0, 0 > \,+\, |~ 1, 1 > ~=~ $ \\
(2.1) \\
$~~~~~~~~~~~~~~~ ~=~ |~ 0 > \bigotimes |~ 0 > \,+\, |~ 1 > \bigotimes |~ 1 > \,\,\in
                                            \mathbb{C}^2 \bigotimes \mathbb{C}^2 $ \\

which is well known to be entangled, in other words, $|~ \omega_{00} >$ is not of the form \\

$~~~~~~ ( \alpha |~ 0 > + \beta |~ 1 > ) \bigotimes( \gamma |~ 0 > + \delta |~ 1 > ) \in
        \mathbb{C}^2 \bigotimes \mathbb{C}^2 $ \\

for any $\alpha, \beta, \gamma, \delta \in \mathbb{C}^2$. \\

Here we can turn to the usual and rather picturesque description used in Quantum Computation,
where two fictitious personages, Alice and Bob, are supposed to exchange information, be it of
classical or quantum type. Alice and Bob can each take their respective qubit from the
entangled pair of qubits $|~ \omega_{00} >$, and then go away with their respective part no
matter how far apart from one another. And the two qubits thus separated in space will remain
entangled, unless of course one or both of them get involved in further classical or quantum
interactions. For clarity, however, we should note that the single qubits which Alice and Bob
take away with them from the pair $|~ \omega_{00} >$ are neither one of the terms $|~ 0, 0 >$
or $|~ 1, 1 >$ above, since both these are themselves already pairs of qubits, thus they
cannot be taken away as mere single qubits, either by Alice, or by Bob. Consequently, the
single qubits which Alice and Bob take away with them cannot be described in any other form,
except that which is implicit in (2.1). \\
Now, after that short detour into the language of Quantum Computation, we can note that,
according to Quantum Mechanics, the entanglement in the double qubit $|~ \omega_{00} >$
implies that the states of the two qubits which compose it are correlated, no matter how far
from one another Alice and Bob would be with them. Consequently, knowing the state of one of
these two qubits can give information about the state of the other qubit. On the other hand,
in view of General, or even Special Relativity, such a knowledge, say by Alice, cannot be
communicated to Bob faster than the velocity of light. \\
And yet, anybody who is familiar enough with Quantum Mechanics, can instantly know and
understand all of that, no matter how far away from one another Alice and Bob may be with
their respective single but entangled qubits. \\

So that, again, the question arises : how and where does such a thinking happen ? \\ \\

{\bf 3. Making Reasonable the Unreasonable ... ?} \\

As seen above, res cogitans has the outmost stubborn nature not to let itself be forgotten,
omitted, disregarded, and let alone dismissed by any sort of ... never mind mind ... \\

And then, perhaps, we may as well try to reconsider our present approaches to the much
puzzling and remarkably powerful fact of our ability to think, and think especially within the
scientific realms. After all, the point made by Wigner, [4], is still of major, and also ever
growing concern, to the extent that technology is so fundamental in modern life, and with it
so is Mathematics, with its ever more unsettling ... unreasonable effectiveness in the Natural
Sciences ... \\
And beyond all utilitarian reasons, the issue of the deeper sources of our thinking is a major
issue of the very human condition, since as far as we know it at present, we humans are the
only creatures capable of such a mode of being, a mode that is one of the most important of
our distinguishing features. \\

In this regard, a Descartes had it relatively easy, being able to see the two distinct
branches of res cogitans and res extensa as being firmly, vitally and eternally sourced in the
idea of a unique, all underlying and unifying Creator, an idea so much prevalent at his
time. \\

In our own days, with these two branches long severed by now from any underlying unity, and
thus standing in an ever unsolvable dichotomy, how could we try to bring them together again
into a deeper unifying source, one that would not necessarily require a mere simple return to
ideas and visions of past times ? \\

Recent developments may suggest a way in this regard, [1]. Indeed, since ancient times, the
self-referential ability of human awareness or consciousness is not only known, but also seen
as just about the most important of its features, [5]. \\
However, and rather amusingly, this self-referentiality has so far, and throughout the ages,
only appeared to be an instant and major source of trouble in a variety of human ventures, be
they usual language, Logic, Mathematics itself, and so on. For instance, the ancient Greek
riddle of the liar, as much a Russell's paradox, or G\"{o}del's proof of his celebrated
theorem on incompleteness are such instances. \\
The conclusion which thus has emerged and got evermore strengthened was to avoid so called
"circular reasoning", a reasoning which contains any sort of self-reference ... \\

And yet, recently, a rigorous mathematical development precisely of such circular or
self-referring thinking has appeared, [1]. \\

Consequently, it may perhaps be the time to try to bring res cogitans, that realm so stubbornly
refusing to be dismissed, together with res extensa which, nowadays, so much occupies the
center stage of our all too utilitarian concerns and ventures. \\

Is this, then, a much better, and long missing approach to mathematical knowledge when seen in
its so called internal dimensions ? \\
Time may of course tell, to the extent that such an avenue would eventually be pursued ... \\

Anyhow, the endless and evermore detailed, varied, learned and sophisticated pursuit of what a
Voegelin found only to be a curiosity shop would not seem much more capable of bringing us
nearer to a more reasonable view of that rather dramatic statement of Wigner on the
unreasonable effectiveness of Mathematics in Natural Sciences ... \\


\begin{thebibliography}{99}

\bibitem{} Barwise J, Moss L : Vicious Circles, on the Mathematics of Non-Wellfounded
Phenomena. CSLI Lecture Notes No. 60, Center for the Study Language and Information,
Stanford, California, 1996

\bibitem{} Manin Yu I : Mathematical knowledge : internal, social and cultural aspects.
arXiv:math.HO/0703427

\bibitem{} Rosinger E E : Where and how does it happen ? \\ arXiv:physics/0505041

\bibitem{} Wigner E : The unreasonable effectiveness of Mathematics in the Natural Sciences.
Comm. Pure and Appl. Math., Vol, 13, No. 1, February 1960

\bibitem{} Genesis 3:14

\end{thebibliography}
\end{document}